\documentclass[ejs]{imsart}

\arxiv{math.PR/0000000}

\usepackage{amsfonts}
\usepackage{amssymb}
\usepackage{graphicx}
\usepackage{amsmath}
\usepackage[a4paper]{geometry}
\usepackage{color}

\newtheorem{theorem}{Theorem}

\newtheorem{remark}{Remark}

\definecolor{jens}{rgb}{0,0.4,0.8}

\begin{document}

\begin{frontmatter}
\title{On signal detection and confidence sets for low rank inference problems}
\runtitle{Low rank inference}

\begin{aug}
\author{\fnms{Alexandra} \snm{Carpentier}\ead[label=e1]{carpentier@maths.uni-potsdam.de}\thanksref{t1}}
\thankstext{t1}{This work was carried out when this author was a research associate in the University of Cambridge.}
\address{Institut f\"ur Mathematik\\
Universit\"at Potsdam\\
Am Neuen Palais 10\\
14469 Potsdam\\
Germany\\
\printead{e1}}
\affiliation{Universit\"at Potsdam}

\and
\author{ \fnms{Richard} \snm{Nickl}\ead[label=e2]{r.nickl@statslab.cam.ac.uk}}
\address{Statistical Laboratory\\
Center for Mathematical Sciences\\
University of Cambridge\\
Wilberforce Road\\
CB3 0WB Cambridge\\
United Kingdom\\
\printead{e2}}
\affiliation{University of Cambridge}

%\thankstext{t1}{Some comment}
\runauthor{A. Carpentier, R. Nickl}

%\phantom{E-mail:\ }\printead*{e2}}

\end{aug}

\begin{abstract}
We consider the signal detection problem in the Gaussian design trace regression model with low rank alternative hypotheses. We derive the precise (Ingster-type) detection boundary for the Frobenius and the nuclear norm. We then apply these results to show that honest confidence sets for the unknown matrix parameter that adapt to all low rank sub-models in nuclear norm do not exist. This shows that recently obtained positive results in \cite{CGEN15} for confidence sets in low rank recovery problems are essentially optimal.

\end{abstract}



%\tableofcontents

%\begin{keyword}[class=AMS]
%\kwd[Primary ]{60K35}
%\kwd{60K35}
%\kwd[; secondary ]{60K35}
%\end{keyword}

%\begin{keyword}
%\kwd{sample}
%\kwd{\LaTeXe}
%\end{keyword}

\end{frontmatter}

\section{Introduction}

Consider the Gaussian design trace regression model
\begin{equation} \label{model}
Y_i=\mathrm{tr}(X^i\theta) + \epsilon_i,~~~~i=1, \dots, n,
\end{equation}
where $\epsilon \sim N(0, I_n)$ is an i.i.d.~vector of Gaussian noise. Here the matrices $X^i$ are $d\times d$ square matrices with i.i.d.~entries $X^i_{mk} \sim N(0,1)$, and $\theta$ is the unknown $d \times d$ matrix we want to make inference on.  We are interested in the case where the model dimension $d^2$ is possibly large compared to sample size $n$, but where $\theta$ has low rank $k$, in which case we write $\theta \in R(k), 1 \le k \le d$. This setting serves as a prototype for various matrix inference problems such as those occurring in compressed sensing \cite{CP11} or in quantum tomography \cite{GLFBE10}. We consider here a high-dimensional regime where $\min(d,n) \to \infty$, reflecting contemporary statistical challenges. 

\medskip

The first problem we study in this paper is the \textit{signal detection problem} with low-rank alternatives: We want to test the hypothesis $$H_0: \theta = 0 ~~vs.~~H_1: \theta \neq 0, \theta \in R(k), \|\theta\| \ge \rho,$$ where $\|\cdot\|$ equals either the Frobenius norm $\|\cdot\|_F$ or the nuclear norm $\|\cdot\|_*$ (defined in detail below), and where $\rho$ should be the minimal `signal strength' condition for the above hypothesis testing problem to have a consistent solution (in the sense of Ingster, see \cite{I03}). We will show that the minimax optimal detection boundary in Frobenius norm is of the form $$\rho \approx \min \left(\sqrt{\frac{d}{n}}, n^{-1/4}\right)$$ whereas in nuclear norm it is $$\rho \approx  \min \left(\sqrt{\frac{kd}{n}}, \sqrt {\frac{k}{n^{1/2}}}\right).$$ A remarkable feature is that for the Frobenius norm the detection rate \textit{does not depend at all} on the complexity of the alternative hypothesis (the rank $k$), whereas for the nuclear norm it does. The phase transition between the two regimes in these rates depends precisely on whether the sample size $n$ exceeds the dimension $d^2$ of the maximal parameter space $R(d)$ or not. The upper bounds in our proofs are related to the papers \cite{ITV10, ACCP11} about the detection boundary in the sparse regression setting, and our main contribution consists in deriving the matching lower bounds for low rank alternatives.

\medskip

Our interest in the detection boundary is triggered by the second problem we investigate here: the question of existence and non-existence of adaptive confidence sets for low rank parameters. It follows from general decision-theoretic principles (see Chapter 8.3 in \cite{GN15} and also \cite{HN11, BN13}) that the answer to this question is closely related to a `composite version' of the detection problem (see (\ref{compo}) below). This approach was employed in \cite{NvdG13} to prove that adaptive and honest confidence sets for the parameter $\theta$ do not exist in \textit{sparse} regression models if an $\ell_2$-risk performance beyond $O(n^{-1/4})$ is desired. In contrast in the recent paper \cite{CGEN15} it was shown that if sparsity constraints are replaced by low rank conditions, then adaptive and fully honest confidence sets exist over the entire parameter space $R(d)$. Adaptation means here that the expected Frobenius norm diameter of the confidence set reflects the minimax risk over arbitrary low rank sub-models $R(k), 1 \le k \le d$. The fact that the detection rates obtained here in Frobenius norm are independent of the rank constraint $\theta \in R(k)$ provides another heuristic explanation of the result in \cite{CGEN15}.

\smallskip

Moreover \cite{CGEN15} constructed another confidence set whose diameter adapts to low rank sub-models in the stronger nuclear norm distance, and that is honest for all $\theta$'s that are \textit{non-negative definite and have trace equal to one}, that is, whenever $\theta$ is the density matrix of a \textit{quantum state}. Such a constraint on $\theta$ is natural in a quantum physics context considered in \cite{CGEN15}, but not in general. The question arises whether it is essentially necessary or not. In the present paper we show that indeed the existence results of \cite{CGEN15} are specific to the geometry induced by the Frobenius norm or to the quantum state constraint, and that nuclear-norm adaptive and honest confidence sets over general low rank parameter spaces \textit{do not exist} in the model (\ref{model}). For example, our results imply that if one requires coverage of a confidence set over all of $R(d)$ then the worst case nuclear norm diameter for rank-one parameters can be off the minimax estimation rate over $R(1)$ by as much as $\sqrt d$. Our results thus further illustrate the subtleties involved in the theory of confidence sets for high-dimensional parameters, and that the positive results in \cite{CGEN15} are of a rather specific nature.

\smallskip

Our proofs are given in the simplest model where both the design and the noise are Gaussian, and the matrices involved are of square type. As usual, our results extend without major difficulty to sub-Gaussian design and noise, to certain correlated random designs, and also to non-square matrices, at the expense of slightly more technical proofs. Generalisations of our results to the matrix completion problem are currently under investigation.

\section{Main results}

\subsection{Notation}
We write $\mathbb M_d$ for the set of $d \times d$ matrices with real elements. If $\mathcal X: \mathbb M_d \to \mathbb R^n$ denotes the `sampling operator' $$\theta \mapsto \mathcal X\theta = \big(\mathrm{tr}(X^1\theta), \dots, \mathrm{tr}(X^n\theta) \big)^T,$$ then the model (\ref{model}) can be written as 
\begin{equation*}
Y = \mathcal X\theta + \epsilon,
\end{equation*}
where $Y=(Y_1, \dots, Y_n)^T$ and $\epsilon = (\epsilon_1, \ldots, \epsilon_n)^T$. We write $E^X$ for the expectation over the distribution of $\mathcal X$ only, and $E_\theta$ for the expectation conditional on $\mathcal X$. The full expectation is denoted by $\mathbb E_\theta= E^X E_\theta$. The corresponding probability laws are denoted by $P^X, P_\theta, \mathbb P_\theta$ and we employ the usual $o/O/o_P/O_P$-notation with $\min(n,d) \to \infty$.

We denote the standard norm on Euclidean space by $\|\cdot\|_2$, and the associated inner product by $\langle \cdot, \cdot \rangle_2$. Let $\|.\|_F$ be the Frobenius norm over $\mathbb M_d$, i.e.
$$\|M\|_F = \sqrt{\mathrm{tr}( M^TM)} = \sqrt{\sum_{j\leq d} \lambda_j^2},$$
where $\lambda_i^2$ are the eigenvalues of $ M^TM$. The associated inner product is $$\langle U, V \rangle_F = \mathrm{tr}(U^TV).$$ We also define the nuclear norm of $M$ as
$$\|M\|_{*} = \sum_{j\leq d} |\lambda_j|.$$ These two norms are in fact defined also for matrices that are not of square type. Finally we recall that for any matrix $M \in R(k)$, we have $$\|M\|_{F}\leq \|M\|_{*}\leq \sqrt{k}\|M\|_{F}.$$

\subsection{Signal detection for low rank alternatives}

We consider first the following hypothesis testing problem, also known as the signal detection problem: 
\begin{equation} \label{maintest}
H_0: \theta = 0~~vs.~~H_1: \theta \in R(k), \|\theta\| \ge \rho.
\end{equation}
Here the alternative space is restricted to a `low rank' hypothesis $\theta \in R(k)$ for some $1 \le k \le d$. Moreover, for a separation constant $\rho>0$, the detection boundary is described by a `signal strength' condition measured in terms of the size $\|\theta\| \ge \rho$ of the Frobenius-, or of the nuclear norm of $\theta$. In the high-dimensional regime where $\min (n,d) \to \infty$, we want to find the minimal sequence $\rho \equiv \rho_{n,d}$ such that for any $\alpha>0$ a level $\alpha$-test $\Psi=\Psi(Y,\mathcal X, \alpha)$ exists: 
\begin{equation} \label{power}
 \left[\mathbb E_0 [\Psi] + \sup_{\theta \in H_1} \mathbb E_\theta [1-\Psi] \right] = \mathbb P_0 (\text{reject } H_0) + \sup_{\theta \in H_1} \mathbb P_\theta (\text{accept } H_0) \le \alpha.
\end{equation}
Recall that a test is simply a random indicator function $\psi=1_A$ where the rejection event $A$ depends only on $Y, \mathcal X, \alpha$, and we require the sum of the type-one and the type-two error of the test to be controlled at any fixed level $\alpha>0$.

\begin{theorem}\label{signalthm} Consider the testing problem (\ref{maintest}) with norm $\|\cdot\|$. Define

\begin{equation*}
 r_{n,d} = \begin{cases}
\min(\sqrt{d/n}, n^{-1/4})&\text{if $\|\cdot\|=\|\cdot\|_F$}\\
 \min(\sqrt{kd/n},\sqrt k/n^{1/4})&\text{if $\|\cdot\|=\|\cdot\|_*$.}
  \end{cases}
\end{equation*}

\medskip

1) Suppose $\rho \ge Dr_{n,d}$. Then for every $\alpha>0$ there exists a test $\Psi= \Psi(Y,\mathcal X, \alpha)$ and finite constants $D=D_\alpha>0, n_\alpha \in \mathbb N$ such that (\ref{power}) holds for every $n \ge n_\alpha$.

\medskip

2) Conversely, suppose $\rho = o(r_{n,d})$ and $k=o(d)$ as $\min(n,d) \to \infty$. Then no test satisfying (\ref{power}) for every $\alpha>0$ exists. In fact 
\begin{equation}\label{ing}
\liminf_{n,d} \inf_{\Psi} \left[\mathbb E_0 [\Psi] + \sup_{\theta \in H_1} \mathbb E_\theta [1-\Psi] \right] \ge 1
\end{equation}
where the infimum extends over all test functions $\Psi= \Psi(Y,\mathcal X)$.

\end{theorem}

\smallskip

The tests $\Psi$ constructed in the proof are given in (\ref{test}) below and straightforward to implement. Note also that the $\|\cdot\|_*$-separated alternatives are a subset of the $\|\cdot\|_F$-separated alternatives (see (\ref{nucsep}) below), and our results imply that an optimal test for the case $\|\cdot\|=\|\cdot\|_F$ is essentially optimal also for $\|\cdot\|_*$.

\subsection{Confidence sets for low rank recovery}

Low rank recovery algorithms are well-studied in compressed sensing and high-dimensional statistics, see e.g., \cite{CP11, GLFBE10, K11, KLT11, NW11, CZ15} and the references therein. In the setting of model (\ref{model}) they provide minimax optimal estimators $\tilde \theta$ of $\theta \in R(k)$ with (high probability) performance guarantees  
\begin{equation} \label{minimax}
\|\tilde \theta - \theta\|^2_F \lesssim \frac{kd}{n},~~~\|\tilde \theta- \theta\|_* \lesssim k \sqrt{\frac{d}{n}}.
\end{equation}
The question we study here is whether associated uncertainty quantification methodology exists, that is, whether we can find confidence sets $C_n \subset \mathbb M_d$ such that 
\begin{equation}\label{coverage}
\inf_{\theta \in \mathbb M_d}\mathbb P_\theta (\theta \in C_n) \ge 1-\alpha,
\end{equation}
 at least for $\min(n,d)$ large enough, and such that the diameter $|C_n|$ of $C_n$ reflects the accuracy of adaptive estimation in the sense that $|C_n|$ shrinks, with high probability, at the optimal rates from (\ref{minimax}) whenever $\theta \in R(k)$. We insist here on an \textit{adaptive} confidence set that does not require knowledge of the unknown rank $k$ of $\theta$.

A first result that is proved in the paper~\cite{CGEN15} is that such adaptive confidence sets do exist in the model (\ref{model}) if the diameter is measured in Frobenius distance. The construction of this set is straightforward, see \cite{CGEN15} for details.

\begin{theorem}[Theorem 2 in \cite{CGEN15}]
For every $\alpha>0$ there exists a confidence set $C_n=C_n(Y,\mathcal X, \alpha)$ such that for all $n \in \mathbb N$, (\ref{coverage}) holds, and such that uniformly in $\theta \in R(k_0)$ for any $1 \le k_0 \le k,$ with high $\mathbb P_\theta$-probability the Frobenius-norm diameter $|C_n|_F$ of $C_n$ satisfies $$|C_n|_{F} \lesssim \sqrt{k_0\frac{d}{n}}.$$  
\end{theorem}

A second result that  is proved in the paper~\cite{CGEN15} is that an (asymptotic) adaptive confidence set exists also in nuclear norm provided that the ``quantum state constraint" is satisfied, namely, provided it is known a priori that $\theta$ is non-negative definite and has nuclear norm one, and provided the coverage requirement in (\ref{coverage}) is relaxed to hold only over a maximal model $R(k)$ in which asymptotically consistent estimation of $\theta$ is possible (i.e., $k\sqrt{d/n} = o(1)$). Define $$R^+(k) = R(k) \cap \{\theta \text{ is non-negative definite},~ \text{tr}(\theta)=1\},$$ the set of quantum state density matrices of rank at most $k$.

\begin{theorem}[Theorem 4 in \cite{CGEN15}]\label{thm:rich}
Assume $k\sqrt{d/n} = o(1)$ for some $1 \le k \le d,$ and let $\alpha>0$ be given. Then there exists a confidence set $C_n=C_n(Y,\mathcal X, \alpha)$ such that
$$\liminf_{\min (n,d) \to \infty} \inf_{\theta \in R^+(k)} \mathbb P_\theta (\theta \in C_n) \ge 1-\alpha,$$ 
and such that uniformly in $\theta \in R^+(k_0)$ for any $1 \le k_0 \le k,$ with high $\mathbb P_\theta$-probability the nuclear norm diameter $|C_n|_*$ of $C_n$ satisfies $$|C_n|_{*} \lesssim k_0\sqrt{\frac{d}{n}}.$$
\end{theorem}

In fact it is not difficult to generalise the above theorem to the case where the condition $tr(\theta)=1$ is relaxed to $\|\theta\|_* \le 1$.

\bigskip

The next theorem, which is the main result of this subsection, implies that no analogue of Theorem 2 can hold true if the Frobenius norm there is replaced by the nuclear norm, and it also shows that Theorem 3 cannot hold true if $R^+(k)$ is replaced by $R(k)$, that is, if the `quantum state constraint' is relaxed. More precisely, we show that if a confidence set $C_n$ is required to have coverage over the maximal model $R(k_1)$, then the worst case expected nuclear norm diameter of $C_n$ over arbitrary sub-models $R(k_0), k_0=o(k_1),$ depends on the maximal model dimension $k_1$ and does not improve as $k_0 \downarrow 1$. The proof of Theorem \ref{main} is based on Part 2) of Theorem \ref{signalthm} and lower bound techniques for adaptive confidence sets from \cite{HN11, BN13}.

\begin{theorem} \label{main}
Let $k_1 \to \infty$ such that $k_1=o(d)$ as $\min(n,d)\to \infty$. Suppose that for any $0<\alpha<1/3$ the confidence set $C_n=C_n(Y, \mathcal X, \alpha)$ is asymptotically honest over the maximal model $R(k_1)$, that is, it satisfies 
\begin{equation}\label{cov2}
\liminf_{\min (n,d) \to \infty} \inf_{\theta \in R(k_1)}\mathbb P_\theta (\theta \in C_n) \ge 1-\alpha.
\end{equation}
Then for every $k_0 =o(k_1)$ and some constant $c>0$ depending on $\alpha$, we have 
\begin{equation} \label{adap}
\sup_{\theta \in R(k_0)} \mathbb E_\theta |C_n|_* \ge c \sqrt{\frac{k_1 d}{n}}
\end{equation}
for every $\min(n,d)$ large enough. In particular no confidence set exists that is honest over all of $\mathbb M_d$ and that adapts in nuclear norm to any model $R(k_0), k_0 =o(\sqrt d)$.  
\end{theorem}

For notational simplicity we have lower bounded the \textit{expected diameter} $|C_n|_*$ in (\ref{adap}), but the proof actually contains a stronger `in probability version' of this lower bound.

\smallskip

\begin{remark}  A few remarks on Theorem \ref{main} are in order: 

\smallskip

i) In the least favourable case where one wants coverage over the entire $R(d) = \mathbb M_d$ while still adapting to rank-one matrices (i.e., $k_0=1$), the performance of any honest confidence set is off the minimax optimal adaptive estimation rate $\sqrt{d/n}$ over $R(1)$ by a diverging factor that can be as close to $\sqrt d$ as desired.  

\smallskip

ii) Even if one restricts coverage to hold only for `consistently estimable models' $R(k_1)$ with $k_1 \sqrt {d/n}\to 0$ (as in Theorem \ref{thm:rich}), the diameter $|C_n|_*$ can be off the minimax rate of estimation over $R(1)$ by a factor of $\sqrt k_1$. 

\smallskip

iii) We also note that the above result does \textit{not} disprove the existence of adaptive confidence sets for sub-models $R(k_0)$ of `moderate rank' where $k_0 \ge \sqrt d$. While more of technical interest -- note that this rules out $n < d^2$ for consistent recovery to be possible -- this regime currently remains open (it is related to the apparently hard problem of finding optimal separation rates in the composite testing problem (\ref{compo}) below).
\end{remark}

\section{Proofs}

\subsection{Proof of Theorem \ref{signalthm}, upper bounds}

When $n<d^2$ then define $$\hat r_n = \frac{1}{n}\|Y\|_2^2 - 1,~~\tau_n = n^{-1/2}$$ but when $n \ge d^2$ set $$\hat r_n = \frac{2}{n(n-1)} \sum_{i<j} \sum_{1 \le m \le d, 1 \le k \le d}Y_i X^i_{mk}Y_jX^j_{mk},~~\tau_n=d/n.$$ The test statistic is 
\begin{equation}\label{test}
\Psi_n = 1\left\{\hat r_n \ge z_\alpha \tau_n \right\}
\end{equation}
where $z_\alpha$ are quantile constants chosen below.

These tests work for Frobenius norm separation, by effectively the same proofs as in \cite{ITV10}, using that we can embed the matrix regression model into a vector regression model with $p=d^2$ parameters, and since the separation rates only depend on the model dimension (and not on low rank or sparsity degrees). However, to provide intuition, we give some details, first for the case $n<d^2$: Under $H_0$ we have $Y=\epsilon$ and so
$$\mathbb E_0 \Psi_n = \Pr\left(\frac{1}{\sqrt n} \sum_{i=1}^n (\varepsilon_i^2 -E \varepsilon_i^2) >z_\alpha \right) \le \alpha/2$$ for every $n \in \mathbb N$ and $z_\alpha$ large enough (using either Chebyshev's inequality and $E\varepsilon_i^4 =3$, or Theorem 4.1.9 in \cite{GN15} for a more precise non-asymptotic bound). Now for the alternatives $\theta \in H_1$ we use the basic concentration result Lemma 1a) in \cite{CGEN15} which implies that for any fixed $\theta$ the event $$\mathcal E = \left\{\left|(1/n)\|\mathcal X \theta\|_2^2 - \|\theta\|_F^2 \right|\le \|\theta\|_F^2/2 \right\}$$ has $P^X$-probability at least $1-2\exp(-n/24)$, and so, for $n\ge n_\alpha$ such that $2\exp(-n/24) <\alpha/6$,
\begin{align*}
\mathbb E_\theta(1-\Psi_n) &= \mathbb P_\theta \left(\hat r_n < z_\alpha \tau_n  \right) \\
&= \Pr \left(\frac{1}{n}\|\mathcal X \theta + \epsilon\|_2^2 -1 < \frac{z_\alpha}{\sqrt n} \right) \\
&= \Pr \left(\frac{1}{n}\|\mathcal X\theta\|_2^2 - \frac{z_\alpha}{\sqrt n} < -\frac{2}{n} \epsilon^T\mathcal X\theta - \frac{1}{n}\sum_{i=1}^n( \varepsilon_i^2-1) \right) \\
&\le \Pr \left(\frac{\|\theta\|_F^2}{2} - \frac{z_\alpha}{\sqrt n} < -\frac{2}{n} \epsilon^T\mathcal X\theta - \frac{1}{n}\sum_{i=1}^n( \varepsilon_i^2-1) , \mathcal E\right) + 2 \exp(-n/24)  \\
&\le \Pr \left(\left|\frac{2}{n} \epsilon^T\mathcal X\theta\right| > \|\theta\|_F^2/8, \mathcal E\right) + \Pr\left(\frac{1}{\sqrt n} \sum_{i=1}^n (\varepsilon_i^2 -E \varepsilon_i^2) >z_{\alpha/3} \right) +\alpha/6
\end{align*}
since, by the hypothesis on $\rho$, we have for $D$ large enough that $$\frac{\|\theta\|_F^2}{2} - \frac{z_\alpha}{\sqrt n} \ge \frac{\|\theta\|_F^2}{4} \ge \frac{2z_{\alpha/3}}{ n^{1/2}}.$$ The last probability is bounded by $\alpha/6$ as under $H_0$ and the last but one probability is also bounded by $\alpha/6$ by a direct (conditional on $\mathcal X$) Gaussian tail inequality (restricting to the event $\mathcal E$: just as in term II of the proof of Theorem \cite{CGEN15} with $\tilde \theta=0$ there), so that in total we have bounded the testing errors in (\ref{power}) by $\alpha/2 + (3/6)\alpha =\alpha$, as desired. The case $n \ge d^2$ follows from similar but slightly more technical arguments, adapting the arguments from proof of Theorem 3 in \cite{CGEN15}, or arguing directly as in Theorem 4.3 in \cite{ITV10} with $p=d^2$.

\medskip

The test (\ref{test}) also works for nuclear-norm separation since $$H_1^* = \theta \in R(k): \|\theta\|_* \ge c\sqrt{k} \rho$$ is a subset of $$H_1^F = \theta \in R(k) : \|\theta\|_F \ge c \rho$$ in view of the inequality 
\begin{equation} \label{nucsep}
\|\theta\|_F \ge (1/\sqrt k) \|\theta\|_* ~~\forall \theta \in R(k),
\end{equation}
 so that $$\mathbb E_0\Psi_n + \sup_{\theta \in H_1^*} \mathbb E_\theta(1-\Psi_n) \le \mathbb E_0\Psi_n + \sup_{\theta \in H_1^F} \mathbb E_\theta(1-\Psi_n) \le \alpha.$$ We now turn to the more difficult lower bounds.

\subsection{Proof of Theorem \ref{signalthm}, lower bounds}

Let $\Psi$ be any test -- any measurable function of $Y, \mathcal X$ that takes values in $\{0,1\}$.  Assume $\rho=o(r_{n,d})$ as $\min(n,d) \to \infty$ and let $H_1=H_1(\rho)$ be the corresponding alternative hypothesis. 

\medskip

\textit{Step I: Reduction to averaged likelihood ratios}: Let $\pi=\pi_{n,d}$ be a sequence of finitely supported probability distributions on $\mathbb M_d$ such that $\pi_{n,d}(H_1) \to 1$, and denote by $\pi|H_1$ that measure restricted to $H_1$ and re-normalised to unit mass. Define $$Z= \mathbb E_{\theta \sim \pi} \prod_{i\leq n} \frac{dP_{i}^{(\theta)}}{dP_{i}^{(0)}} \equiv \int \prod_{i\leq n} \frac{dP_{i}^{(\theta)}}{dP_{i}^{(0)}} d\pi(\theta) ,$$ where $dP_{i}^{(\theta)}$ is the distribution of $Y_i|\mathcal X$ when the parameter generating the data is $\theta$, and $dP_{i}^{(0)}$ is the distribution of $Y_i|\mathcal X$ when the parameter generating the data is $0$. Then, by a standard testing lower bound (e.g., (6.23) in \cite{GN15}), for any $\eta>0$,
\begin{align*}
\mathbb E_0 \Psi + \sup_{\theta \in H_1} \mathbb E_\theta(1-\Psi)   &\ge \mathbb E_0 \Psi + \mathbb E_{\theta \sim \pi|H_1} \mathbb E_\theta(1-\Psi)\\
& \ge \mathbb E_0 \Psi + \mathbb E_{\theta \sim \pi}\mathbb E_\theta(1-\Psi) - o(1) \\
&= E^X\left[E_0 \Psi + \mathbb E_{\theta \sim \pi}  E_\theta(1-\Psi) \right] - o(1) \\
&\geq (1-\eta)\left[1 - \left[\frac{\sqrt{\mathbb E_0(Z-1)^2}}{\eta}\right]\right] - o(1). 
\end{align*}
Now since $$\mathbb E_0[Z-1]^2 = \mathbb E_0[Z^2] -1,$$ if we show that $\mathbb E_0[Z^2]\le 1+o(1)$ as $\min(n,d) \to\infty$ for a suitable choice of $\pi$, then the lower bound (\ref{ing}) will follow by letting $\eta \to 0$. Recall the notation $\mathbb E_\theta = E^X E_\theta$.

\medskip

\textit{Step II: Computation of $E_0[Z^2]$:} The $(Y_i)$ are independent with distribution $\mathcal N((\mathcal X\theta)_i,1)$ conditional on the design $\mathcal X$, hence
\begin{align*}
Z &=\mathbb E_{\theta \sim  \pi} \Bigg[ \prod_{i \leq n} \frac{\exp(-\frac{1}{2}(y_{i} - (\mathcal X\theta)_i)^2)}{\exp(-\frac{1}{2}y_{i}^2)} \Bigg]\\
&=\mathbb E_{\theta \sim \pi}\Bigg[ \prod_{i \leq n} \exp(y_{i} (\mathcal X\theta)_i) \exp(- \frac{1}{2} ((\mathcal X\theta)_i)^2) \Bigg]
\end{align*}
and can hence write
\begin{small}
\begin{align}
 E_0 \big[Z^2 \big] &= \int_{\mathbb R^n} \Bigg(\mathbb E_{\theta \sim \pi} \Big[\prod_{i \leq n} \exp(y_{i}  (\mathcal X\theta)_i) \exp(- \frac{1}{2} ((\mathcal X\theta)_i)^2)\Big] \Bigg)^2 \prod_{i\leq n} \frac{1}{\sqrt{2\pi}}\exp(- \frac{y_i^2}{2} ) dy_1...dy_{n} \nonumber\\
&= \int_{\mathbb R^n} \Bigg(\mathbb E_{\theta \sim \pi} \Big[ \exp(- \frac{1}{2} \|\mathcal X \theta\|_2^2) \prod_{i \leq n} \exp(y_{i}  (\mathcal X \theta)_i \Big] \Bigg)^2 \prod_{i\leq n} \frac{1}{\sqrt{2\pi}}\exp(- \frac{y_i^2}{2} ) dy_1...dy_{n} .\nonumber
\end{align}
\end{small}
Thus, if $\theta, \theta'$ are independent copies of joint law $\pi^2$, then we have
\begin{align*}
E_0 \big[Z^2 \big]  &= \int_{\mathbb R^n} \mathbb E_{\pi^2} \Big[  \exp(- \frac{1}{2} (\|\mathcal X \theta\|_2^2 - \frac{1}{2} (\|\mathcal X \theta'\|_2^2)\prod_{i \leq n} \frac{1}{\sqrt{2\pi}}\exp\big(y_{i} (\mathcal X(\theta+\theta'))_i - \frac{y_i^2}{2}\big) \Big] dy_1...dy_{n}\\
&= \mathbb E_{\pi^2} \Bigg[ \exp(- \frac{1}{2} \|\mathcal X \theta\|_2^2 - \frac{1}{2} \|\mathcal X\theta'\|_2^2)\nonumber\\  
&~~~~\times \prod_{i \leq n}\int_{y_i} \Big(\frac{1}{\sqrt{2\pi}} \exp\Big(-\frac{1}{2} \big(y_{i} - (\mathcal X(\theta+ \theta'))_i)^2\Big)dy_i \exp\Big(\frac{1}{2} (\mathcal X(\theta+\theta'))_i^2\Big)\Bigg]\nonumber\\
&= \mathbb E_{\pi^2}\left[ \exp\left(\frac{1}{2}\|\mathcal X(\theta + \theta')\|_2^2 - \frac{1}{2} \|\mathcal X \theta\|_2^2 - \frac{1}{2} \|\mathcal X \theta'\|_2^2 \right)\right]\nonumber
%&= \mathbb E_{\theta, \theta'}\Bigg[ \exp\Big(\langle \mathcal X(\theta), \mathcal X(\theta')\rangle\Big)\Bigg].\label{coli}
\end{align*}

\medskip

\textit{Step III: Integrating over $\mathcal X$:} The $E^X$-expectation of the last expression can be bounded by 
$$
\mathbb E_{\pi^2} \left[\exp \left(\frac{n}{2} (\|\theta+\theta'\|_F^2 - \|\theta\|_F^2 - \|\theta'\|_F^2) \right) E^X \exp \left(\frac{1}{2}(Z_1-Z_2-Z_3) \right)\right]
$$
where $$Z_\ell = \|\mathcal X\vartheta_\ell\|_2^2 - n \|\vartheta_\ell\|_F^2,~~\text{with}~~\vartheta_1=\theta+\theta', \vartheta_2=\theta, \vartheta_3=\theta'.$$ The last factor can be bounded, by applying the Cauchy-Schwarz inequality twice, by
\begin{equation} \label{product}
(E^X\exp(Z_1))^{1/2} (E^X \exp(2Z_2))^{1/4} (E^X \exp(2Z_3))^{1/4}.
\end{equation}
Since $\mathcal X \vartheta_\ell \sim N(0, \|\vartheta_\ell\|_F^2 I_n)$ the distribution of $Z_\ell$ is the one of $\|\vartheta_\ell\|_F^2 \sum_{i=1}^n (g_i^2-1)$ where the $g_i$ are i.i.d.~$N(0,1)$. Applying Theorem 3.1.9 in \cite{GN15} with $\tau_i \equiv 1$ and $\lambda = \|\vartheta_1\|_F^2$ or $\lambda = 2\|\vartheta_\ell\|_F^2, \ell=2,3,$ (and hence setting $\|A\|=1, \|A\|_{HS}=n$ in that theorem) we see that if $\max_\ell \|\vartheta_\ell\|_F^2 \le 1/4$ then
$$E^X \exp(Z_1) \le \exp\left(\frac{n \|\vartheta_1\|_F^4}{1- 2 \|\vartheta_1\|_F^2} \right),~~\text{and}~~E^X \exp(2Z_\ell) \le \exp\left(\frac{2n \|\vartheta_\ell\|_F^4}{1- 4 \|\vartheta_\ell\|_F^2} \right),~\ell=2,3.$$ As a consequence if 
\begin{equation} \label{viertel}
\max_{\ell=1,2,3} \|\vartheta_\ell\|_F = o(n^{-1/4})
\end{equation} 
then the the product (\ref{product}) is bounded above by $1+o(1)$. We conclude that if the prior $\pi$ satisfies (\ref{viertel}) almost surely then
\begin{align*}
\mathbb E_0[Z^2] = E^XE_0[Z^2] &\le \left(1+o(1) \right) \times E_{\pi^2} \exp \left(\frac{n}{2} (\|\theta+\theta'\|_F^2 - \|\theta\|_F^2 - \|\theta'\|_F^2) \right)  \\
& = \left(1+o(1) \right) \times E_{\pi^2} \exp \left(n \langle \theta, \theta' \rangle_F \right) .
\end{align*}

\medskip

\textit{Step IV: Construction of $\pi$ and bounds for $\mathbb E_0[Z^2]$:} Assume for notational simplicity that $d$ is an integer multiple of $k$, the general case needs only minor notational adjustment. Pick independent random $d \times 1$ vectors $v_\ell: \ell=1, \dots, k$ each of which consists of i.i.d.~Rademacher entries (i.e., taking values $\pm 1$ with probability $1/2$). Create a matrix $W$ as follows: In the first $d/k$ columns insert $v_1$ times a random sign $B_{1,j}, j=1, \dots, d/k$. Then, in the $\ell$-th block repeat the same with $v_1$ replaced by $v_\ell$, and random signs $B_{\ell,j}, j=1, \dots, d/k$. If $\|\cdot\|=\|\cdot\|_F$ let $\gamma_n=\rho_n/d $ and if $\|\cdot\|=\|\cdot\|_*$ set $\gamma_n = 2\rho_n /(\sqrt k d)$, so that in either case $$\gamma_n= o\left(\min(\sqrt{1/dn}, d^{-1}n^{-1/4}) \right).$$  Define the random matrix $\theta = \gamma_n W$ and let $\theta'$ be an independent copy of it. Thus $$n\langle \theta, \theta' \rangle_F = n\gamma_n^2 \sum_{\ell =1}^k \sum_{m=1}^d \sum_{j=1}^{d/k} v_{\ell,m} B_{\ell,j} v'_{\ell,m} B_{\ell,j}' = n\gamma_n^2 \sum_\ell \sum_m v_{\ell,m} v_{\ell,m}' \sum_j B_{\ell,j} B'_{\ell,j}.$$  As products of Rademacher variables are again Rademacher variables we have, for $\epsilon_{\ell, m},\tilde \epsilon_{\ell,j}$ i.i.d.~Rademacher variables (all defined on a suitable product probability space),
\begin{align} \label{tb}
E_{\pi^2} \exp \left(n \langle \theta, \theta' \rangle_F \right) &= E_\epsilon E_{\tilde \epsilon} \exp \left(n\gamma_n^2 \sum_{\ell} \sum_{m} \epsilon_{\ell, m} \sum_j \tilde \epsilon_{\ell, j} \right) \notag \\
&= \left(E_\epsilon E_{\tilde \epsilon} \exp \left(n\gamma_n^2 \sum_{m} \epsilon_{\ell, m} \sum_j \tilde \epsilon_{\ell, j} \right)\right)^k.
\end{align} Conditional on the values of $\epsilon$ we set $\lambda = n \gamma_n^2 \sum_{m=1}^d \epsilon_{\ell,m}$ and note that $$|\lambda| \le nd \gamma_n^2 = o(1).$$ By Taylor expansion or standard properties of the hyperbolic cosine (as, e.g., in the proof of Theorem 6.2.9 in \cite{GN15})
$$E_{\tilde \epsilon} \exp\left(\lambda \sum_{j=1}^{d/k} \tilde \epsilon_{\ell, j}\right) = \cosh (\lambda^2)^{d/k} \le \exp\left(\lambda^2 d/k  \right)$$ and thus, since $[EU]^k \le E[U^k]$ for any non-negative random variable $U$, the right hand side in (\ref{tb}) is bounded above by
$$ \left(E_\epsilon  \exp\left(\lambda^2 d/k \right) \right)^k \le E_\epsilon \exp\left(\lambda^2 d\right) =E_\epsilon \exp\left(n^2 \gamma_n^4 d \left(\sum_{m=1}^d \epsilon_m\right)^2 \right) \equiv  E \exp \left(Z^2/c^2 \right)$$ where the Rademacher sum $Z=\sum_{m=1} ^d \epsilon_m$ is a sub-Gaussian random variable with variance proxy $\sigma^2=d$ (cf.~Section 2.3 in \cite{GN15}). Thus by (2.24) in \cite{GN15} we have
$$E \exp \left(Z^2/c^2 \right) \le 1 + \frac{2}{c^2/2\sigma^2 -1} = 1+o(1)$$
since
$$\frac{c^2}{\sigma^2} = \frac{1}{d^2 n^2 \gamma_n^4} \to \infty$$ as $n,d \to \infty$. Summarising all steps so far we conclude $$0 \le E^XE_0[Z-1]^2 = \mathbb E[Z^2]  -1 \le 1-1+o(1)=o(1)$$ noting that (\ref{viertel}) holds $\pi$-almost surely in view of $$\|\theta\|_F^2 = \gamma_n^2 \|W\|_F^2 = \gamma_n^2 d^2 = o(n^{-1/2}).$$

\medskip

\textit{Step V: Asymptotic concentration of $\pi$ on $H_1$}: Finally we show that for the above prior we have indeed $\Pi(H_1)\to 1$. First since $\theta$ consists of columns that are linear combinations of at most $k$ distinct vectors $v_\ell$ we immediately have $\theta \in R(k)$ almost surely. Moreover, for the case $\|\cdot\|=\|\cdot\|_F$ we have from the last display and by definition of $\gamma_n$ that $\|\theta\|_F^2 = \rho^2_n$, so $\Pi(H_1)=1$ follows. 

For the case $\|\cdot\|=\|\cdot\|_*$ we have to show that $$\pi_{n,d}(\|\theta\|_* \ge \rho_n) \to 1$$ as $\min(n,d) \to \infty$. We can transform $\theta$ into the $d \times k$ matrix $\theta U$ consisting of $k$ column vectors $\gamma_n \sqrt{d/k} v_\ell, \ell=1, \dots, k$. The corresponding $d \times k$ matrix $U$ consists of $k$ column vectors, the $\ell$-th of which has zero entries except for the indices $m \in [\ell d/k, \dots, -1+(\ell+1)d/k]$, where it equals $\sqrt {k/d} B_{\ell, m}$. Thus, $U$ is an orthonormal projection matrix and we deduce that $$\|\theta\|_* \ge \|\theta U\|_*.$$ We can renormalise the column vectors of $\theta U$ so that $$\theta U = \gamma_n \frac{d}{\sqrt k} \left(\dots \frac{1}{\sqrt d}v_\ell \dots \right) \equiv \gamma_n \frac{d}{\sqrt k} V.$$ The $d \times k$ matrix $V$ consists of scaled i.i.d.~Rademacher entries, and hence the proof of Lemma 1 in \cite{NvdG13} (with $n=d, k=k_1=p$ in the first display on p.2868 there) implies that, if $k/d \to 0$, then with probability as close to one as desired, the smallest singular value of $V$ is bounded below by $1/2$ for $d$ large enough. As a consequence $\|V\|_* \ge k/2$ and so, with probability approaching one, $$\|\theta\|_* \ge \gamma_n d\sqrt{k}/2 =\rho_n.$$ Note that the same lower bound holds for 
\begin{equation}\label{compdist}
\|\theta-R(k_0)\|_* =\inf_{\theta' \in R(k_0)} \|\theta-\theta'\|_* \ge \sum_{j=k_0+1}^k|\lambda_j| \ge (k-k_0)/2
\end{equation} 
for any $k_0<k$, if the absolute eigenvalues in the last display are assumed to be in decreasing order.

\subsection{Proof of Theorem \ref{main}}

Consider the composite testing problem 
\begin{equation}\label{compo}
H_0: \theta \in R(k_0) ~vs ~H^c_1: \theta \in R(k_1), \|\theta - R(k_0)\|_*=\inf_{\theta' \in R(k_0)} \|\theta-\theta'\|_* \ge \rho.
\end{equation} 
From (\ref{compdist}) with $k=k_1$ and $k_0 = o(k_1)$ we see that for $\min(n,d)$ large enough such that $(k_1-k_0)/2 \ge k_1/4$, the prior $\pi$ from the previous proof with $\gamma_n=4\rho_n /(\sqrt k d)$ asymptotically concentrates on $H_1^c$. As a consequence testing (\ref{compo}) is no easier than when $H_0=\{0\}$, so that when $\rho=o(\sqrt{k_1 d/n})$ then the proof of Part 2 of Theorem \ref{signalthm} implies 
\begin{equation}\label{errors}
\liminf_{n,d} \inf_{\Psi} \left[\sup_{\theta \in H_0} \mathbb E_\theta\Psi + \sup_{\theta \in H_1^c} \mathbb E_\theta(1-\Psi)\right] \ge 1.
\end{equation} 
Now assume by way of contradiction that there exists $C_n$ that satisfies (\ref{cov2}) with $\alpha<1/3$ and such that for every $c>0$ there exist infinitely many $n,d$ such that $$\sup_{\theta \in H_0}\mathbb E_\theta |C_n|_*<c \sqrt{k_1d/n}.$$ Passing to the infinite subsequence $\min(n,d) \to \infty$ along which the last inequalities hold, we deduce from Markov's inequality that $$\sup_{\theta \in R(k_0)}\mathbb P_\theta(|C_n|_* \ge \alpha \sqrt{k_1d/n}) \le c/\alpha<\alpha$$ for $c$ small enough depending only on $\alpha$. Then, by Proposition 8.6.3 in \cite{GN15} we can construct a test for (\ref{compo}) for which the testing errors in (\ref{errors}) are no more than $3 \alpha<1$ along the chosen subsequence, a contradiction that completes the proof.

\medskip

\textbf{Acknowledgement.} RN's research was supported by the European Research Council (ERC) under grant agreement No.647812. This paper was written while AC was a research associate in the University of Cambridge. The authors are grateful to two referees whose comments improved the exposition of this article. 

%\bibliographystyle{plain}
%\bibliography{reference3}

\end{document}